\documentclass[12pt,a4paper]{article}
\usepackage[cp1251]{inputenc}
\usepackage[ukrainian]{babel}
\usepackage[cp1251]{inputenc} % Кодова сторінка 1251.
\usepackage[ukrainian]{babel} % Використовуємо українську та англійську (за замовченням) мови.
\usepackage{bezier,graphics,graphicx}
\usepackage{amsmath,amsfonts,amssymb,amscd}
\usepackage[mathcal]{eucal}
\usepackage{latexsym}
\usepackage{layout}
\usepackage{makeidx}
\usepackage{epsfig}

\usepackage{amsmath,amsthm,amscd,amssymb,color}
\usepackage{latexsym}
\usepackage{epsf}
\usepackage{epsfig}
\usepackage{floatflt}

\newcommand{\R}{\mathbb{R}}
\newcommand{\rank}{\operatorname{rank}}
\newcommand{\N}{\mathbb{N}}
\newcommand{\pd}{\partial}

\theoremstyle{plain}
\newtheorem{theorem}{Theorem}[section]
\newtheorem{lemma}[theorem]{Lemma}

\newtheorem{proposition}[theorem]{Proposition}

\theoremstyle{definition}

\newtheorem{example}[theorem]{Example}
\theoremstyle{remark}

\newtheorem{case[theorem]}{Case}

\begin{document}
\large
\noindent
УДК 517.9; 519.6\\  MSC 34K28, 65N22, 15A03, 34B05

\hfill {Математика}

\bigskip

\noindent
%\begin{center}
{\bf  Oksana Bihun${}^{\ast}$ and Mykola Prytula${}^{\ast\ast}$ \footnote{\copyright Oksana Bihun and Mykola Prytula} }
%\end{center}

\vskip8pt

\noindent
%\begin{center}

\noindent {\Large \bf  Rank of Projection-Algebraic Representations of Some\\ Differential Operators}

\noindent

\vskip8pt

{\it Abstract.}
The Lie-algebraic method approximates differential operators that are formal polynomials of $\{1,x,\frac{d}{dx}\}$ by linear operators acting on a finite dimensional space of polynomials. In this paper we prove that the rank of the $n$-dimensional representation of the operator $$K=a_k \frac{d^k}{dx^k}+a_{k+1}\frac{d^{k+1}}{dx^{k+1}}+\ldots +a_{k+p}\frac{d^{k+p}}{dx^{k+p}}$$ is $n-k$ and conclude that the Lie-algebraic reductions of differential equations allow to approximate only \emph{some} of solutions of the differential equation $K[u]=f$. We show how to circumvent this obstacle when solving boundary value problems by making an appropriate change of variables. We generalize our results to the case of several dimensions and illustrate them with numerical tests.

\vskip20pt

\section{Introduction}

The Lie-algebraic method of discrete approximations was proposed by Calogero in his seminal work 
\cite{CalogeroInterpNC},  with the purpose of solving eigenvalue problems for ordinary and partial differential operators.  Several years later, Mytropolski, Prykarpatsky, and Samoylenko have established a general algebraic-projection framework for the method that has broadened its applicability to a wide class of problems of Mathematical Physics~\cite{MPSUMZh}. 
 The fundamental idea of the method is to project functions onto the space $P_n[x]$ of polynomials of fixed degree $n$ and to approximate differential operators by linear maps acting on $P_n[x]$. In early 2000s, Bihun and Lu\'{s}tyk estimated the error of such approximation and applied the method to solving boundary value problems for elliptic partial differential equations~\cite{BMatSt, BLVisnyk, BLMatSt}. An overview of the literature on the Lie-algebraic method of discrete approximations can be found in~\cite{BPNTSh}.

In this paper we prove that the rank of the $n$-dimensional representation of the operator $$K=a_k \frac{d^k}{dx^k}+a_{k+1}\frac{d^{k+1}}{dx^{k+1}}+\ldots +a_{k+p}\frac{d^{k+p}}{dx^{k+p}}$$ is $n-k$ and conclude that the Lie-algebraic reductions of differential equations allow to approximate only \emph{some} of solutions of the differential equation $K[u]=f$. We show how to circumvent this obstacle when solving boundary value problems by making an appropriate change of variables. We generalize our results to the case of several dimensions and illustrate them with numerical tests.

\section{The Lie-Algebraic Method\\ of Discrete Approximations}

Every differentiable function $u:[a,b]\to\R$ can be approximated by its Lagrange polynomial of degree $n$. To construct the Lagrange polynomial $L[u]$ of a differentiable function $u:[a,b]\to\R$, the interval $[a,b]$ is partitioned by points $a=x_0<x_1<\ldots<x_n=b$ and the values $L[u](x_i)=u(x_i)$ are prescribed, where
 $i=0,1,2,\ldots,n$. The Lagrange polynomial of $u$ is uniquely defined by its $n+1$ values at the points of the partition and can be written in the form
 $$
 L[u](x)=\sum_{i=0}^n u(x_i)l_i(x),
 $$
 where $$l_i(x)=\frac{\prod_{\stackrel{k=0}{k\neq i}}^n (x-x_k)}{\prod_{\stackrel{k=0}{k\neq i}}^n (x_i-x_k)}$$ is the $i$-th standard Lagrange polynomial, $0\leq i\leq n$.
 
 Let $A$ be a differential operator acting on differentiable functions defined on $[a,b]$.  We assume that $A$ belongs to the enveloping algebra of the Heisenberg-Weyl algebra of differential operators $\{1,x,\frac{d}{dx}\}$; that is, $A$ is a formal polynomial of $\{1,x,\frac{d}{dx}\}$. The value of $A$ at a function $f$ is denoted by $A[f]$. For a differentiable function $u:[a,b] \to \R$ and a partition $\sigma=\{x_0,x_1,\ldots,x_n\}$ of $[a,b]$, let $u_\sigma=(u(x_0),u(x_1),\ldots,u(x_n))^T$ be the vector of the values of $u$ at the nodes of the partition. The finite dimensional representation $A_\sigma:\R^{n+1}\to\R^{n+1}$ of $A$ can be constructed from the condition
 $$
A_\sigma u_\sigma=A[L[u]]_\sigma.
 $$

 Following this condition, the finite dimensional representation $Z=\{Z_{kj}\}_{k,j=0}^n$ of the operator $d/dx$ can be found from the relation
$$
\left(
\begin{array}{c}
\frac{dl_k}{dx}(x_0)\\
\frac{dl_k}{dx}(x_1)\\
\ldots\\
\frac{dl_k}{dx}(x_n)
\end{array}
\right)=Z
\left(
\begin{array}{c}
l_k(x_0)\\
l_k(x_1)\\
\ldots\\
l_k(x_n)
\end{array}
\right),
$$
where $0\leq k\leq n$. Because $l_k(x_j)=\delta_{kj}$ for all $k,j \in \{0,1,\ldots,n\}$, the right hand side of the last equality is the $k$-th column of $Z$. Therefore,
\begin{equation}
\label{df:Z}
 Z_{jk}=\frac{dl_k}{dx}(x_j)=
 \left\{
 \begin{array}{ll}
 \sum_{m=0, m\neq j}^n \frac{1}{(x_j-x_m)} \text{ if } k=j,\\
 \frac{\pi_j}{\pi_k}\frac{1}{x_j-x_k} \text{ if } k\neq j,
 \end{array}
 \right.
\end{equation}
where $\pi_k=\prod_{m=0, m\neq k}^n(x_k-x_m)$ for every $k\in\{1,2,\ldots,n\}$.

Because the derivative $\frac{dl_k}{dx}$ coincides with its Lagrange polynomial, we obtain the relation
\begin{equation}
\label{eq:dlkZ}
\frac{dl_k}{dx}(x)=\sum_{j=0}^n Z_{jk} l_j(x).
\end{equation}

The finite dimensional representation $X$ of the operator $x$ of multiplication by the variable $x$ is the square matrix of dimension $(n+1)$ with the diagonal $(x_0,x_1,\ldots,x_n)$:
\begin{equation}
\label{df:X}
X=\operatorname{Diag}(x_0,x_1,\ldots, x_n).
\end{equation}
 The finite dimensional representation of the identity operator $1$ is the $(n+1)$-dimensional identity matrix denoted by $I_{n+1}$.

Having constructed the finite dimensional representations of the operators $1,x,d/dx$, we can approximate every operator $A$ in the enveloping algebra of $\{1,x,d/dx\}$. 

\begin{example}
\label{ex1}
Consider the  boundary value problem
\begin{equation}
\label{eq:ex1}
\left\{
\begin{array}{ll}
u''+u=0,\\
u(0)=2, u(\pi/2)=1
\end{array}
\right.
\end{equation}
for a $C^2$ function $u:[0,\pi/2]\to\R$. Let $\sigma=\{x_i\}_{i=0}^n$ be a partition of $[0,\pi/2]$ and let $Z$ and $I_{n+1}$ be the finite dimensional representations of $\frac{d}{dx}$ and $1$, respectively, for this partition.  We approximate the differential operator $A=d^2/dx^2+1$ by $Z^2+I_{n+1}$ and reduce problem~\eqref{eq:ex1} 
 to the problem of finding a vector $u_\sigma \in \R^{n+1}$ that satisfies 
\begin{equation}
\label{eq:ex1approx}
\left\{
\begin{array}{ll}
(Z^2 +I_{n+1})u_\sigma=0,\\
{}[u_\sigma ]_0=2, [u_\sigma]_n=1.
\end{array}
\right.
\end{equation}
\end{example}

A natural question to ask here is this: does problem~\eqref{eq:ex1approx} have a solution? More generally, when does the finite dimensional approximation of the original problem have a unique solution in $\R^{n+1}$? We address this question in the following sections.

\section{Rank of the Finite Dimensional Representation of the Differential Operator $\sum_{p=k}^m a_p\frac{d^p}{dx^p}$}
\label{sc:Rank1d}
Let $1\leq k\leq m\leq n$ and let $a_k, a_{k+1},\ldots, a_m$ be real numbers such that $a_k$ and $a_m$ are not zeroes. 
In this section we prove that the finite dimensional representation 
$$a_k Z^k+a_{k+1} Z^{k+1}+\ldots +a_m Z^m$$ of the differential operator 
$$a_k \frac{d^k}{dx^k}+a_{k+1} \frac{d^{k+1}}{dx^{k+1}}+\ldots +a_m \frac{d^m}{dx^m}$$ has rank $n+1-k$.

%We denote the set ring of all polynomials of degree $n$ with real coefficients by $P_n[x]$.

\begin{lemma}
\label{lemma1}
Let $\sigma=\{ x_0,x_1,\ldots,x_n \}$ 
be a partition of the interval $[a,b]$ 
and let $Z$, given by formula \eqref{df:Z}, be the finite dimensional representation of $d/dx$ with respect to this partition. Then 
$ \rank Z=n$ and $Z^{n+1}=0$.
\end{lemma}
\textbf{Proof.}
Let us recall that the matrix $Z$ is $(n+1)$-dimensional.

First, let us show that  $Z^{n+1}=0$. Using formula~\eqref{eq:dlkZ}, it is easy to verify that 
$$
\frac{d^k l_j}{dx^k}(x)=\sum_{m=0}^n [Z^k]_{mj} l_m(x)
$$
for every $k \in \N$ and $j \in \{0,1,\ldots,n\}$. Because $\frac{d^{n+1} l_j}{dx^{n+1}}\equiv 0$, we obtain 
$$
\sum_{m=0}^n [Z^{n+1}]_{mj}l_m(x)=0,
$$ for all $x \in [0,\pi/2]$, which, after evaluation $x=x_k$ for $k\in\{0,1,\ldots,\}$, implies that $Z^{n+1}=0$.

Let us now prove that $ \rank Z=n$. 

There exists a nonsingular matrix $B=\{B_{ij}\}_{i,j=0}^n$ of dimension $n+1$ such that 
$$
l_j(x)=\sum_{p=0}^n B_{pj}x^p
$$
for every integer $j$ between $0$ and $n$.
Then 
\begin{eqnarray*}
\frac{dl_j}{dx}(x)&=&\sum_{m=0}^n Z_{mj} l_m(x)\\
&=&\sum_{m=0}^n \sum_{p=0}^n Z_{mj} B_{pm} x^p\\
&=&\sum_{p=0}^n [BZ]_{pj} x^p.
\end{eqnarray*}

Because the polynomial $\frac{dl_j}{dx}(x)$ has degree $n-1$, the last row of the matrix $BZ$ must be zero. Therefore, $\rank Z=\rank BZ\leq n$.

It remains to show that $\rank Z$ cannot be strictly smaller than $n$. Suppose, on the contrary, that $\rank Z<n$. Then there exist $n$ linearly dependent columns in $Z$. Without loss of generality we assume that the first $n$ columns of $Z$ are linearly dependent; denote these columns by $A_0, A_1, \ldots A_{n-1}$. There exist constants $\alpha_0, \alpha_1, \ldots, \alpha_{n-1}$, not all zeros, such that
$$
\sum_{j=0}^{n-1} \alpha_j A_j=0.
$$
But then 
\begin{eqnarray*}
\sum_{j=0}^{n-1} \alpha_j \frac{dl_j}{dx}(x)&=&\sum_{j=0}^{n-1}\sum_{m=0}^n \alpha_j [A_j]_m l_m(x)\\
&=&\sum_{m=0}^n[\sum_{j=0}^{n-1} \alpha_j A_j]_m l_m(x)\\
&=&0,
\end{eqnarray*}
 which, after integration, gives us 
$$\sum_{j=0}^{n-1} \alpha_j l_j(x)=C,$$ where $C$ is a constant.
 
 By setting $x=x_j$ in the previous equality, where $j \in\{0,1,\ldots,n-1\}$, we obtain $\alpha_0=\alpha_1=\ldots=\alpha_{n-1}=C \neq 0$. On the other hand, if we set $x=x_n$ in the very same equality, we obtain $C=0$, a contradiction. We have thus proved that $\rank Z$ cannot be smaller than $n$. Q.E.D.
 
 \begin{lemma}
 \label{lemma2}
 Let $H$ be a nilpotent square matrix of dimension $n+1$ such that $\rank H=n$ and $H^{n+1}=0$. Then $\rank H^k=n+1-k$ for all $k \in\{0,1,\ldots,n\}$. Moreover, the Jordan form of $H$ is 
\begin{equation}
\label{JordH}
J=\Big(\begin{array}{cc}
0 & I_n\\
0& 0
\end{array}\Big),
\end{equation}
where $I_n$ is the identity matrix of dimension $n$.
 \end{lemma} 
\textbf{Proof.}
Let $J$ be the Jordan form of $H$ so that $H=B^{-1} J B$ for some invertible matrix $B$ and $\operatorname{diag}(J)=(\lambda_0, \lambda_1,\ldots,\lambda_n)$, where $\lambda_i$ are the eigenvalues of $H$.

Then $0=H^{n+1}=B^{-1} J^{n+1} B$ and $J^{n+1}=0$. Because $\operatorname{diag}(J^{n+1})=(\lambda_0^{n+1}, \lambda_1^{n+1},\ldots,\lambda_n^{n+1})$, the matrix $J$ has only one eigenvalue $\lambda=0$ of multiplicity $n+1$ and is given by formula~\eqref{JordH}. 

It is now easy to verify that 
\begin{equation}
J^k=\Big(\begin{array}{cc}
0 & I_{n+1-k}\\
0& 0
\end{array}\Big)
\end{equation}
and $\rank J^k=n+1-k$. Thus, $\rank H^k=\rank J^k=n+1-k$ as required. Q.E.D.

\begin{theorem}
\label{theorem1}
Let $P(z)=a_k z^k+ a_{k+1} z^{k+1}+\ldots+a_m z^m$ be a polynomial of degree $m$ with real coefficients, where $k \in\{0,1,\ldots,m\}$ and $a_k\neq 0$. If a square matrix $B$ is nilpotent, then $\rank P(B)=\rank B^k$.
\end{theorem}

\textbf{Proof.}

Let us rewrite $P(B)$ as $P(B)=a_k B^k(I+\frac{a_{k+1}}{a_k} B+\frac{a_{k+2}}{a_k} B^2+\ldots +\frac{a_m}{a_k}B^{m-k})$. Define $D:=\frac{a_{k+1}}{a_k} B+\frac{a_{k+2}}{a_k} B^2+\ldots +\frac{a_m}{a_k}B^{m-k}$. 

Because the matrix $B$ is nilpotent, $D$ is nilpotent. But then the matrix $I+D$ is invertible. Indeed, the series $\sum_{p=0}^\infty (-D)^p$ is a finite sum and equals $(I+D)^{-1}$.

In summary, $P(B)=B^k A$, where $A$ is an invertible matrix. Therefore, $\rank P(B)=\rank B^k$ as required. Q.E.D.

The next theorem follows immediately from theorem~\ref{theorem1} and lemmas~\ref{lemma1} and~\ref{lemma2}.

\begin{theorem}
\label{theorem3}
Let $\sigma=\{ x_0,x_1,\ldots,x_n \}$ 
be a partition of the interval $[a,b]$ 
and let $Z$, given by formula~\eqref{df:Z}, be the finite dimensional representations of  the operators $d/dx$ and $x$, respectively, corresponding to this partition. Let $P(z)$ be a polynomial defined in the hypothesis of theorem~\ref{theorem1}. Then 
$ \rank Z=n$, $Z^{n+1}=0$, and $\rank P(Z)=\rank{Z^k}=n+1-k$.
\end{theorem}

It is natural to ask whether theorem~\ref{theorem1} can be generalized to the case where the coefficients $a_p$, $p\in\{k,k+1,\ldots,m\}$, of the polynomial $P(z)=a_k z^k+ a_{k+1} z^{k+1}+\ldots+a_m z^m$, are polynomials in $x$. In other words, suppose that $a_p$, $p\in\{k,k+1,\ldots,m\}$, are polynomials, matrix $B$ is nilpotent and matrix $X$ is diagonal. Is it true that the matrix $P(X,B)$ has the rank of $B^k$? The next example shows that, in general, this is not true.

\begin{example}
Let $$B=
\left(\begin{array}{cc}
-2 & -1\\
4 & 2
\end{array}\right), \;\mbox{ and }\;\; X=
\left(\begin{array}{cc}
a & 0\\
0 & b
\end{array}\right).$$
It is easy to verify that $B^2=0$, hence $B$ is nilpotent. Consider the polynomial $P(x,z)=1+xz$. The matrix $P(X,B)=I_2+XB$ does not have the full rank if $b-a=-\frac{1}{2}$. Indeed,
$$I_2+XB=
\left(\begin{array}{cc}
1-2a & -a\\
4b & 1+2b
\end{array}\right)
$$ and $\operatorname{det}(I_2+XB)=1+2(b-a)=0$ if $b-a=-\frac{1}{2}$.
\end{example}

Let us revise example~\ref{ex1} in the view of what we have proved. The finite-dimensional approximation of the differential equation $u''+u=0$ on the interval $[0,\pi/2]$ is given by 
\begin{equation}
\label{la1}
Z^2 u_\sigma+u_\sigma=0,
\end{equation} where $u_\sigma\in \R^n$ is a vector whose coordinates approximate the values of the solution $u$ at the nodes of a partition $\sigma$ of $[0,\pi/2]$.

By the last theorem, the rank of the matrix $Z^2+I$ is $n+1$, thus the system of linear equations~\eqref{la1} has a \emph{unique} solution. On the other hand, we know that the space of solutions to the equation $u''+u=0$ is two-dimensional. Therefore, such direct application of the Lie-algebraic method does not help us approximate all the solutions of the differential equation, only some of them. It is possible, however, to still approximate solutions of initial and boundary value problems using the Lie-algebraic method. We simply need to make an appropriate substitution so that to ensure that the solution of the differential equation automatically satisfies initial or boundary conditions.

For example, to solve the boundary value problem
\begin{equation}
\label{BVP2}
\left\{
\begin{array}{ll}
u''+u=0, \; x \in[0,\frac{\pi}{2}]\\
u(0)=2, u(\frac{\pi}{2})=1.
\end{array}
\right.
\end{equation}
we make a substitution $u=(2-2/\pi x)\big(x(x-\pi/2)v+1\big)$ and obtain the equation $p(x) v''+q(x) v'+r(x) v=s(x)$, where
\begin{eqnarray*}
p(x)&=&-\frac{2}{\pi} x^3+3x^2-\pi x,\\
q(x)&=&2(-\frac{6}{\pi} x^2+6x-\pi),\\
r(x)&=&-\frac{12}{\pi} x+6+x(2-\frac{2}{\pi} x)(x-\frac{\pi}{2}),\\
s(x)&=&\frac{2}{\pi} x-2,
\end{eqnarray*}
 whose discrete approximation is
$$\big(p(X) Z^2 +q(X) Z+r(X) \big) v_\sigma=s(X),$$
where $\sigma$ is a partition of $[0,\pi/2]$.

We have chosen the partition of the interval $[0.001,\frac{\pi}{2}]$ with equally distributed $n$ nodes $x_i=0.001+\frac{i}{n}(\frac{\pi}{2}-0.001)$. We compared the values of the vector $v_{\sigma}$ with the values $\big(u(x_0), u(x_1),\ldots u(x_n)\big)$ of the exact solution $u(x)=\sin x+ 2\cos x$ to the boundary value problem~\eqref{BVP2}. The values of the errors $$E(n)=\sum_{k=0}^n \big|[u_\sigma]_k-u(x_k)\big|$$ and $$E_{max}(n)=\max_{k \in \{0,1,\ldots,n\}} \big|[u_\sigma]_k-u(x_k)\big|$$ are given in table~1.

We have compared the results obtained by means of the Lie-algebraic method with the results obtained using the standard shooting method. Following the shooting technique, we have split boundary value problem~\eqref{BVP2} into two initial value problems:
\begin{equation}
\label{ivp1}
\left \{
\begin{array}{l}
w''=-w\\
w(0)=2, w'(0)=0,
\end{array}
\right.
\end{equation}
and
\begin{equation}
\label{ivp2}
\left \{
\begin{array}{l}
v''=-v\\
v(0)=0, v'(0)=1.
\end{array}
\right.
\end{equation}
The solution of problem~\eqref{BVP2} can be found by the formula $$u=w+\frac{1-w(\frac{\pi}{2})}{v(\frac{\pi}{2})}v.$$

We solved the initial value problems~\eqref{ivp1} and~\eqref{ivp2} using finite differences.

The error of the shooting method for problem~\eqref{BVP2} is given in table~1, which demonstrates fast convergence of the Lie-algebraic scheme.

\begin{center}
\begin{tabular}{c|c|c|c|c|c}
\label{tttv1}
Method & $n$ & 4 & 8 & 12 & 16 \\
\hline
Lie-Alg. & $E(n)$ & 2.2788e-04 & 9.5522e-07 & 3.9033e-09 & 1.5205e-11\\
Lie-Alg. & $E_{max}(n)$ & 1.1466e-04 & 2.5575e-07 & 6.8542e-10 & 1.9955e-12\\
Shooting &  $E(n)$ & 2.71e-02 & 1.39e-02 & 9.3e-03 & 7.0e-03\\
Shooting &  $E_{max}(n)$ & 1.1e-02 & 2.7e-03 & 1.2e-03 & 6.7013e-04
\end{tabular}
%\caption{The errors $E(n)$ and $E_{max}(n)$ of the Lie-algebraic scheme and the shooting method for boundary value 
%problem~\eqref{BVP2}, where $n$ is the number of subintervals in the partition of $[0.001,\frac{\pi}{2}]$ for the %Lie-algebraic scheme and $[0,\frac{\pi}{2}]$ for the shooting method.}
\end{center}

\begin{center}
\small Table~1. The errors $E(n)$ and $E_{max}(n)$ of the Lie-algebraic scheme and the shooting method for boundary value 
problem~\eqref{BVP2}, where $n$ is the number of subintervals in the partition of $[0.001,\frac{\pi}{2}]$ for the Lie-algebraic scheme and $[0,\frac{\pi}{2}]$ for the shooting method.
\normalsize
\end{center}

\section{ The Lie-algebraic Approximations\\ in Several Dimensions}

\subsection{Lagrange Polynomials of Several Variables}
 In this section we will generalize the results proved in the previous section for the case of several dimensions. More precisely, we will study the existence and the number of solutions of the Lie-algebraic discrete approximations of partial differential equations.
 
 Consider a differential equation 
 $$
 Au=f,
 $$
 where $\Omega$ is a $d$-dimensional cube $[a_1, b_1]\times [a_2, b_2] \times \ldots \times [a_d, b_d]$ and $f:\Omega \to \R$ belongs to the domain of the differential operator $A$. We denote $x=(x^{1},\ldots, x^{d})$.
 
 We partition each interval $[a_k, b_k]$, where $1\leq k \leq d$, with the points $a_k=x^{k}{}_0<x^{k}{}_1<\ldots < x^{k}{}_{n_k}=b_k$, and denote the basic Lagrange polynomials with respect to this interval by 
$$l^{k}{}_j(x^{k})=\prod_{m=0, m\neq j}^{n_k} (x^k-x^k{}_m)/\prod_{m=0, m\neq j}^{n_k} (x^k{}_j-x^k{}_m),$$ 
where $j=0,1,\ldots, n_k$.
 
 The basic Lagrange polynomials on $\Omega$ are defined by
 $$
 L_{(i_1, i_2, \ldots ,i_d)}(x^{1},x^{2},\ldots, x^{d})=l^{1}{}_{i_1}(x^{1}) l^{2}{}_{i_2}(x^{2})\ldots l^{d}{}_{i_d}(x^{d}).
 $$
 For brevity, we introduce the multi-index notation $(i)=(i_1, i_2, \ldots i_d)$, where $0\leq i_k \leq n_k$ for every $k=1,2,\ldots, d$. The set of all multi-indices is denoted by $\mathcal{I}$ and the nodes of the partition of $\Omega$ are $x_{(i)}:=(x^1{}_{i_1}, x^2{}_{i_2},\ldots, x^d{}_{i_d})$, where $(i) \in \mathcal{I}$. In this notation, the basic Lagrange polynomials on $\Omega$ can be written as $L_{(i)}(x)$, where $(i) \in \mathcal{I}$, and the Lagrange polynomial of a function $u:\Omega \to \R$ is
 $$L(u)=\sum_{(i) \in \mathcal{I}} u(x_{(i)})L_{(i)}(x)$$.
 
 Let $n_\alpha$-dimensional square matrix $\hat{Z}^{(\alpha)}$ denote the representation of the differential operator $\frac{d}{dx^{\alpha}}$ with respect to the partition $\{x^\alpha{}_k\}_{k=0}^{n_\alpha}$ of $[a_\alpha, b_\alpha]$, where $\alpha=1,2,\ldots,d$. It is easy to verify that for $1\leq \alpha\leq d$
 $$
 \frac{\pd}{\pd x^\alpha} L_{(i)}(x)=\sum_{(j) \in \mathcal{I}} W^{(\alpha)}{}_{(j),(i)}L_{(j)}(x),
 $$
 where
$$W^{(\alpha)}{}_{(j),(i)}=W^{(\alpha)}{}_{(j_1,j_2,\ldots,j_d),(i_1,i_2,\ldots,i_d)}
=\delta_{i_1,\ldots,\hat{i_\alpha},\ldots,i_d}^{j_1,\ldots,\hat{j_\alpha},\ldots,j_d}\hat{Z}^{(\alpha)}{}_{j_\alpha,i_\alpha},$$
 $\delta$ is the Kronecker symbol, and a ``hat'' above the index $\hat{i_\alpha}$ means that the index $i_\alpha$ is omitted.

We consider $W^{(\alpha)}$ as a linear map from $V:=\R^{n_1+1}\otimes \R^{n_2+1}\otimes \ldots \otimes \R^{n_d+1}$ to itself, which acts in the following way. Let $\{e^k{}_j\}_{j=0}^{n_k}$ be the standard basis of $\R^{n_k+1}$. Then $$W^{(\alpha)}(e^1{}_{i_1}\otimes e^2{}_{i_2}\otimes \ldots \otimes e^d{}_{i_d})=\sum_{(j)\in \mathcal{I}}W^{(\alpha)}_{(i),(j)} e^1{}_{j_1}\otimes e^2{}_{j_2}\otimes \ldots \otimes e^d{}_{j_d}.$$ 

We identify the map $W^{(\alpha)}$ with a square matrix $\hat{W}^{(\alpha)}$ of dimension $$N:=(n_1+1)\cdot (n_2+1)\cdot \ldots \cdot (n_d+1)$$ that represents the linear map $W^{(\alpha)}$ in the standard basis of $\R^{N}\cong V$ as follows. We enumerate the multi-indices in $\mathcal{I}$ by introducing the following bijection $\ast:\mathcal{I}\to \{1,2,\ldots,N\}$:
\begin{eqnarray*}
\ast(0,0,0,\ldots,0)&=&1,\\
\ast(1,0,0,\ldots,0)&=&2,\\
\ldots&&\\
\ast(n_1,0,0,\ldots,0)&=&n_1+1,\\
\ast(0,1,0,\ldots,0)&=&(n_1+1)+1,\\
\ldots\\
\ast(n_1,1,0,\ldots,0)&=&(n_1+1)+(n_1+1),\\
\ldots\\
\ast(n_1,n_2,0,\ldots,0)&=&n_2\cdot(n_1+1)+(n_1+1),\\
\ast(0,0,1,\ldots,0)&=&1\cdot(n_2+1)(n_1+1)+1,\\
\ast(n_1,n_2,n_3,\ldots,n_d)&=&n_d\cdot(n_{d-1}+1)\ldots (n_1+1)\\
 &&+\ldots+n_3\cdot(n_2+1)(n_1+1)+n_2\cdot(n_1+1)+n_1+1.
\end{eqnarray*}
It is easy to see that 
\begin{eqnarray*}
\ast(i)&=&\ast(i_1,i_2,\ldots,i_d)\\
&=&i_d(n_1+1)(n_2+1)\ldots(n_{d-1}+1)+\ldots+i_2(n_1+1)+i_1+1
\end{eqnarray*}
for all $(i) \in \mathcal{I}$.

Let $\{g_i\}_{i=1}^N$ be the standard basis of $\R^N$. For every multi-index $(i)\in \mathcal{I}$, we identify the basis element $g_{\ast(i)}$ with the basis element $e^1{}_{i_1}\otimes e^2{}_{i_2}\otimes \ldots \otimes e^d{}_{i_d}$ of~$V$. Using this identification, it is now easy to compute the matrix $\hat{W}^{(\alpha)}$ that represents the linear map $W^{(\alpha)}:V \to V$: its $\ast(i)$-th column equals
\begin{eqnarray*}
\hat{W}^{(\alpha)} g_{\ast(i)}&=&
W^{(\alpha)}(e^1{}_{i_1}\otimes e^2{}_{i_2}\otimes \ldots \otimes e^d{}_{i_d})\\
&=&\sum_{(j)\in \mathcal{I}} W^{(\alpha)}_{(i),(j)}(e^1{}_{j_1}\otimes e^2{}_{j_2}\otimes \ldots \otimes e^d{}_{j_d})\\
&=& \sum_{(j)\in \mathcal{I}} W^{(\alpha)}_{(i),(j)} g_{\ast(j)}.
\end{eqnarray*}
Therefore,
\begin{equation}
\label{eq:WhatFormula}
\hat{W}^{(\alpha)}_{\ast(i),\ast(j)}=W^{(\alpha)}_{(i),(j)}=\delta_{i_1,\ldots,\hat{i_\alpha},\ldots,i_d}^{j_1,\ldots,\hat{j_\alpha},\ldots,j_d}\hat{Z}^{(\alpha)}{}_{i_\alpha,j_\alpha}.
\end{equation}
In the following, we identify $W^{(\alpha)}$ and $\hat{W}^{(\alpha)}$.

\begin{proposition}The matrix representation of the differential operator $\frac{\pd}{\pd x^\alpha}$ is given by the $N \times N$-dimensional matrix
\begin{equation}
\label{eq:WAlpha}
\hat{W}^{(\alpha)}=I_{n_1+1}\otimes \ldots \otimes I_{n_{\alpha-1}+1}\otimes \hat{Z}^{(\alpha)} \otimes I_{n_{\alpha+1}+1}\otimes \ldots \otimes I_{n_d+1},
\end{equation}
where $\otimes$ denotes tensor product and $I_{n_k+1}$ is $(n_k+1)$-dimensional identity matrix.
\end{proposition}
We prove the above proposition for the case $d=2$ in the next subsection, where we discuss the structure of matrix representations of $\displaystyle{\frac{\pd}{\pd x^\alpha}}$ for the two-dimensional case. The proof for multi-dimensional case ($d>2$) is similar and is left to the reader. 

It can also be shown that the finite dimensional representation of the operator of multiplication by a variable $x^{\alpha}$ equals
\begin{equation}
\label{eq:XAlpha}
\hat{X}^{(\alpha)}=I_{n_1+1}\otimes \ldots \otimes I_{n_{\alpha-1}+1}\otimes X^{(\alpha)} \otimes I_{n_{\alpha+1}+1}\otimes \ldots \otimes I_{n_d+1},
\end{equation}
where $X^{(\alpha)}=\operatorname{Diag}(x^{\alpha}{}_0, x^{\alpha}{}_1,\ldots, x^{\alpha}{}_{n_\alpha})$.

We end this subsection with stating useful identities for matrix $\hat{W}^{(\alpha)}$.

\begin{lemma}
\label{lemmaWXAlpha}
The matrix
\begin{equation}
\label{eq:WAlphaK}
[\hat{W}^{(\alpha)}]^k=I_{n_1+1}\otimes \ldots \otimes I_{n_{\alpha-1}+1}\otimes [\hat{Z}^{(\alpha)}]^k \otimes I_{n_{\alpha+1}+1}\otimes \ldots \otimes I_{n_d+1}
\end{equation}
represents differential operator $\frac{\pd^k}{\pd (x^\alpha)^k}$ for all $k \in \{1,2,\ldots,n_\alpha\}$, and
\begin{eqnarray}
\label{eq:WAlphaBeta}
\hat{W}^{(\alpha)}\cdot\hat{W}^{(\beta)}&=&I_{n_1+1}\otimes \ldots \otimes I_{n_{\alpha-1}+1}\otimes \hat{Z}^{(\alpha)} \otimes I_{n_{\alpha+1}+1}\otimes \ldots\\
\nonumber
 &&\otimes I_{n_{\beta-1}+1}\otimes \hat{Z}^{(\beta)} \otimes I_{n_{\beta+1}+1}\otimes
\ldots \otimes I_{n_d+1}
\end{eqnarray}
represents differential operator $\frac{\pd^2}{\pd x^\alpha\pd x^\beta}$ for all $\alpha,\beta \in\{1,2,\ldots,d\}$.

Moreover, the matrices $\hat{W}^{(\alpha)}$, where $\alpha \in \{1,2,\ldots,d\}$, commute and are nilpotent.
\end{lemma}

Formulas~\eqref{eq:WAlphaK} and~\eqref{eq:WAlphaBeta} follow from a composition property of tensor products: For every pair $A, B$ of linear maps (between linear spaces) , $(A \otimes B)^2=(A\otimes B)(A \otimes B)=A^2\otimes B^2$. The last statement of the last lemma is an immediate consequence of equations~\eqref{eq:WAlphaK} and~\eqref{eq:WhatFormula}, and lemma~\ref{lemma1}.

\subsection{Matrix Representations of $\frac{\pd}{\pd x^\alpha}$ \\ for the Two-dimensional Case}
\label{subs:W1dim}

Let us observe the structure of the matrices $\hat{W}^{(1)}$ and $\hat{W}^{(2)}$ in the two-dimensional case ($d=2$):
\begin{equation}
\hat{W}^{(1)}=
\left(
\begin{array}{cccc}
\hat{Z}^{(1)}&0& \ldots & 0\\
0&\hat{Z}^{(1)} & \ldots & 0\\
\ldots & {}& {}&{} \\
0 & 0 & \ldots & \hat{Z}^{(1)}
\end{array}
\right),
\end{equation}
where the number of blocks $\hat{Z}^{(1)}$ on the diagonal of $\hat{W}^{(1)}$ is $n_2+1$;
\begin{equation}
\hat{W}^{(2)}=
\left(
\begin{array}{cccc}
\hat{Z}^{(2)}_{0,0}I_{n_1+1}&\hat{Z}^{(2)}_{0,1}I_{n_1+1}& \ldots & \hat{Z}^{(2)}_{0,n_2}I_{n_1+1}\\
\hat{Z}^{(2)}_{1,0}I_{n_1+1}&\hat{Z}^{(2)}_{1,1} I_{n_1+1}& \ldots & \hat{Z}^{(2)}_{1,n_2}I_{n_1+1}\\
\ldots & {}& {}&{} \\
\hat{Z}^{(2)}_{n_2,0} I_{n_1+1}& \hat{Z}^{(2)}_{n_2,1} I_{n_1+1}& \ldots & \hat{Z}^{(2)}_{n_2,n_2}I_{n_1+1}
\end{array}
\right),
\end{equation}
where $I_{n_1+1}$ is the $(n_1+1)$-dimensional identity matrix.

Let us recall that for two linear operators $A:X \to X$ and $B:Y \to Y$ acting on linear spaces $X$ and $Y$, their tensor product $A \otimes B:X \otimes Y \to X \otimes Y$ is defined by $A \otimes B(x\otimes y):=Ax \otimes By$. 

Let $\{e_0, \ldots, e_{n_1}\}$ and $\{f_0, \ldots, f_{n_2}\}$ be the standard bases of $\R^{n_1+1}$ and $\R^{n_2+1}$, respectively.  Let $N=(n_1+1)(n_2+1)$ and $(i)=(i_1,i_2) \in \mathcal{I}$. Having identified the $\ast(i)$-th vector in the standard basis $\{g_1,g_2,\ldots,g_N\}$ of $\R^{N}$ with $e_{i_1}\otimes f_{i_2}$, it is now easy to observe that $\hat{W}^{(1)}=\hat{Z}^{(1)}\otimes I_{n_2+1}$ and $\hat{W}^{(2)}=I_{n_1+1}\otimes \hat{Z}^{(2)}$.

Indeed, the $\ast(i)$-th column of the matrix $\hat{Z}^{(1)}\otimes I_{n_2}$ equals
\begin{eqnarray*}
[\hat{Z}^{(1)}\otimes I_{n_2+1}]_{\cdot,\ast(i)}&=&
\hat{Z}^{(1)}\otimes I_{n_2+1}(e_{i_1}\otimes f_{i_2})\\
&=&\hat{Z}^{(1)}e_{i_1}\otimes I_{n_2+1}f_{i_2}\\
&=&\sum_{j_1=0}^{n_1} \hat{Z}^{(1)}_{j_1,i_1} e_{j_1}\otimes f_{i_2}\\
&=&\sum_{j_1=0}^{n_1}\sum_{j_2=0}^{n_2} \hat{Z}^{(1)}_{j_1,i_1} \delta_{i_2}^{j_2} e_{j_1}\otimes f_{j_2}\\
&=&\sum_{(j) \in \mathcal{I}} W^{(1)}_{(j),(i)} e_{j_1}\otimes f_{j_2}\\
&=& \sum_{(j) \in \mathcal{I}} \hat{W}^{(1)}_{\ast(j),\ast(i)} g_{\ast(j)},
\end{eqnarray*}
which is the $\ast(i)$-th column of $\hat{W}^{(1)}$. The equality $\hat{W}^{(2)}=I_{n_1+1}\otimes \hat{Z}^{(2)}$ can be verified in a similar fashion.

\subsection{Rank of Matrix Representations of Operators $\frac{\pd}{\pd x^\alpha}$ \\ and Their Compositions}

In this subsection we compute rank of matrix representations of operators $\frac{\pd}{\pd x^\alpha}$, and their compositions, by generalizing results of section~\ref{sc:Rank1d}.

We start with proving the following simple fact.

\begin{lemma}
Let $A:\R^m \to \R^m$ and $B:\R^p \to \R^p$ be linear maps with ranks $m_1$ and $m_2$ respectively. Then $\rank (A \otimes B)=m_1 \cdot m_2$.
\end{lemma}
\textbf{Proof.}
Denote
$$
k_A:=\dim(\ker(A))=m-m_1 \mbox{ and } k_B:=\dim(\ker(B))=p-m_2.
$$
Let $\{f_1, f_2,\ldots,f_{k_A}\}$ be a basis of $\ker(A)$; complete it to a basis $$\{f_1, f_2,\ldots,f_{k_A}, \hat{f}_{k_A+1},\ldots,\hat{f}_{m}\}$$ of $\R^m$. Similarly, let $\{g_1, g_2,\ldots,g_{k_B}\}$ be a basis of $\ker(B)$; complete it to a basis $\{g_1, g_2,\ldots,g_{k_B}, \hat{g}_{k_B+1},\ldots,\hat{g}_{p}\}$ of $\R^p$. For every $u \in \R^m$ and $v \in \R^p$, $A\otimes B(u\otimes v)=Au \otimes Bv=0$ if and only if $Au=0$ or $Bv=0$. Therefore,

\begin{eqnarray*}
\ker (A\otimes B)&=&\operatorname{span}\{f_1\otimes g_i, f_2\otimes g_i,\ldots,f_{k_A}\otimes g_i: i=1,2,\ldots,k_B\}\\
&&+\operatorname{span}\{f_1\otimes \hat{g}_i, f_2\otimes \hat{g}_i,\ldots,f_{k_A}\otimes \hat{g}_i: i=k_B+1,\ldots,p \}\\
&&+\operatorname{span}\{\hat{f}_i\otimes g_1, \hat{f}_i\otimes g_2,\ldots,\hat{f}_i\otimes g_{k_B}: i=k_A+1,\ldots,m \},
\end{eqnarray*}
and 
$$
\dim\big(\ker (A\otimes B)\big)=k_A\cdot k_B+k_A\cdot m_2+k_B \cdot m_1=pm-m_1 m_2,\\
$$
which implies $\rank(A\otimes B)=m_1 m_2$ as required.
Q.E.D.

Using the previous lemma, it is easy to verify the statement of the following proposition.

\begin{proposition} Let  $B$ and $C$ be $(n_\alpha+1)$- and $(n_\beta+1)$-dimensional square matrices, respectively.
Let
$$ 
M:=I_{n_1+1}\otimes \ldots \otimes I_{n_{\alpha-1}+1} \otimes B \otimes
I_{n_{\alpha+1}+1}\otimes \ldots \otimes I_{n_d+1}
$$ and
$$
K:=I_{n_1+1}\otimes \ldots \otimes I_{n_{\alpha-1}+1} \otimes B \otimes
I_{n_{\alpha+1}+1}\otimes \ldots \otimes I_{n_{\beta-1}+1} \otimes C \otimes
I_{n_{\beta+1}+1}\otimes \ldots\otimes I_{n_d+1}.
$$
Then $$\rank M^k=\rank B^k\cdot \frac{(n_1+1)\ldots (n_d+1)}{(n_\alpha+1)}$$ and 
$$\rank K=\rank B \cdot \rank C \cdot \frac{(n_1+1)\ldots (n_d+1)}{(n_\alpha+1)(n_\beta+1)}$$.
\end{proposition}

From the previous proposition, formulas~\eqref{eq:WAlpha},~\eqref{eq:WAlphaK},~\eqref{eq:WAlphaBeta}, and 
lemmas~\ref{lemma1},~\ref{lemma2}, we conclude the following.

\begin{proposition}
Let $\alpha, \beta \in \{1,2,\ldots,d\}$. For all $1\leq k \leq n_\alpha$ and 
$1\leq l \leq n_\beta$, the rank of the matrix representation $[\hat{W}^{(\alpha)}]^k$ of $\frac{\pd^k}{\pd (x^\alpha)^k}$ is 
$$\rank[\hat{W}^{(\alpha)}]^k=(n_\alpha+1-k)\frac{(n_1+1)(n_2+1)\ldots (n_d+1)}{(n_\alpha+1)}$$
and the rank of the matrix representation $[\hat{W}^{(\alpha)}]^k [\hat{W}^{(\beta)}]^l$ of $\frac{\pd^{(k+l)}}{\pd (x^\alpha)^k \pd (x^\beta)^l}$ is 
$$\rank [\hat{W}^{(\alpha)}]^k [\hat{W}^{(\beta)}]^l
=(n_\alpha+1-k)(n_\beta+1-l)\frac{(n_1+1)(n_2+1)\ldots (n_d+1)}{(n_\alpha+1)(n_\beta+1)}.$$
Moreover, the matrix $\hat{W}^{(\alpha)}$ is nilpotent: $[\hat{W}^{(\alpha)}]^{n_\alpha+1}=0$.
\end{proposition}

Let $\sigma=\{x_{(i)}\}_{i \in \mathcal{I}}$ be a partition of $\Omega=[a_1,b_1]\times \ldots \times [a_d,b_d]$ as before. Denote the number of nodes in the partition by $N$. Recall that the matrices $\hat{W}^{(\alpha)}$ and $\hat{X}^{(\alpha)}$ denote the matrix representations of the differential operator $\frac{\pd}{\pd x^{\alpha}}$ and the operator of multiplication by $x^\alpha$, respectively, with respect to the partition $\sigma$ (see equations~\eqref{eq:WAlpha} 
and~\eqref{eq:XAlpha}).
Consider partial differential equation 
\begin{equation}
\label{eq:PDE1}
P(x^1,\ldots,x^d,\frac{\pd}{\pd x^1},  \ldots, \frac{\pd}{\pd x^d})u=f,
\end{equation} 
where $P(x^1,\ldots,x^d,\frac{\pd}{\pd x^1},  \ldots, \frac{\pd}{\pd x^d})$ is a formal polynomial of its arguments. Let $\hat{F}=(f(x_{(\gamma_1)}), f(x_{(\gamma_2)}), \ldots, f(x_{(\gamma_N)}) )$, where $\gamma_k \in \mathcal{I}$ and $\ast(\gamma_k)=k$ for all $k \in \{1,2,\ldots,N\}$.

To solve PDE~\eqref{eq:PDE1} using Calogero method,  we need to solve the system linear equations 
\begin{equation}
\label{eq:SLEq1}
P(\hat{X}^{(1)},\ldots,\hat{X}^{(d)},\hat{W}^{(1)},  \ldots, \hat{W}^{(d)})\hat{U}=\hat{F}
\end{equation} 
for $\hat{U} \in \R^N$; the $i$-th component  $\hat{U}_i$ of the solution vector $\hat{U}$, if it exists, approximates the value $u(x_{(\gamma)})$, where $\ast(\gamma)=i$, of solution $u$ of equation~\eqref{eq:PDE1} at the node $x_{(\gamma)}$ of the partition $\sigma$. It is therefore important to know the conditions under which the system of linear equations~\eqref{eq:SLEq1} has a unique solution.

\begin{theorem}
\label{thm:MDRank}
Let $P(z^1,z^2,\ldots,z^d)$ be a polynomial with real coefficients. The matrix $P(\hat{W}^{(1)},\hat{W}^{(2)},\ldots,\hat{W}^{(d)})$ has full rank if and only if the constant term of $P$ is nonzero.
\end{theorem}
\textbf{Proof.}
(i) If the constant term of $P$, which we denote by $a$, is nonzero, then $P(z^1,z^2,\ldots,z^d)=a+P_1(z^1,z^2,\ldots,z^d)$, where $P_1$ is a polynomial in $(z_1,\ldots,z_d)$ whose constant term vanishes. Because matrices $\hat{W}^{(\alpha)}$, where $\alpha \in\{1,2,\ldots,d\}$, commute and are nilpotent (see lemma~\ref{lemmaWXAlpha}), the matrix $P_1(\hat{W}^{(1)},\ldots,\hat{W}^{(d)})$ is nilpotent. But then the matrix $I_N +\frac{1}{a} P_1(\hat{W}^{(1)},\ldots,\hat{W}^{(d)})$ is nonsingular, as the formal series expansion of $(I_N +\frac{1}{a} P_1(\hat{W}^{(1)},\ldots,\hat{W}^{(d)}))^{-1}$ is a finite sum. Therefore, the matrix $$P(\hat{W}^{(1)},\ldots,\hat{W}^{(d)})=a(I_N +\frac{1}{a} P_1(\hat{W}^{(1)},\ldots,\hat{W}^{(d)}))$$ has full rank as required.

On the other hand, if the constant term of $P$ vanishes, then the matrix $P(\hat{W}^{(1)},\ldots,\hat{W}^{(d)})$ is nilpotent, and therefore does not have full rank. Q.E.D.

\begin{example} We finish this article with an application of the Lie-algebraic method to the following partial differential equation of hyperbolic type~\cite{pma6}
\begin{equation}
\label{eq:lala}
u_{xx}-u_{yy}+y u_x=4(y^2-x^2)\sin (1-x^2-y^2)-2xy\cos (1-x^2-y^2)
\end{equation}
on region $\Omega=\{(x,y): x^2+y^2\leq 1\}$ with the boundary condition
\begin{equation}
\label{eq:lalaBC}
u|_{\pd \Omega}=0.
\end{equation}
This problem has solution $u(x,y)=\sin (1-x^2-y^2)$. We compute the approximate solution of the problem using Lie-algebraic approximations and compare our results with the exact solution.

To incorporate boundary conditions, we make the change of variables $$u=(1-x^2-y^2)v$$, which yields the partial differential equation
\begin{eqnarray}
\label{Ex2eq1}
&&(1-x^2-y^2)(v_{xx}-v_{yy}+y v_x)-4x v_x+4y v_y-2xy v\\
\nonumber
&=&4(y^2-x^2)\sin (1-x^2-y^2)-2xy \cos (1-x^2-y^2).
\end{eqnarray}

We partition the square $D=[-1,1]\times[-1,1]$ that contains region $\Omega$ with points $\sigma=\{ (x_i,y_j): 0\leq i\leq n_1, 0\leq j\leq n_2\}$, where $\sigma_x=\{x_i\}_{i=0}^{n_1}$ and $\sigma_y=\{y_i\}_{i=0}^{n_2}$ are partitions of the interval $[-1,1]$.

Let $Z_x$ and $Z_y$ be matrix representations of $d/dx$ on $[-1,1]$ with respect to partitions $\sigma_x$ and $\sigma_y$, respectively. Similarly, let $X$ and $Y$ be matrix representations of the operator of multiplication by $x$ on $[-1,1]$ with respect to partitions $\sigma_x$ and $\sigma_y$, respectively. Matrices $Z_x$ and $X$ have dimension $n_1+1$, while matrices $Z_y$ and $Y$ have dimension $n_2+1$. Finally, the identity matrices $I_{n_1+1}$ and $I_{n_2+1}$ represent the identity operator with respect to partitions $\sigma_x$ and $\sigma_y$, respectively.

The following table gives matrix representations of the operators $\pd/\pd x$, $\pd/\pd y$, $x$, $y$, and $1$, with respect to partition $\sigma$.

  \begin{center}
\begin{tabular}{l|l|l|l|l|l}
\label{ttt2}
Operator & $\displaystyle{\frac{\pd}{\pd x}}$ &  $\displaystyle{\frac{\pd}{\pd y}}$ & $x$ & $y$ & $1$ \\
\hline
\hline
Repres. &$Z_x\otimes I_{n_2+1}$ &    $I_{n_1+1}\otimes Z_y$ & $X \otimes I_{n_2+1}$ & $I_{n_1+1} \otimes Y$ & $I_{n_1+1} \otimes I_{n_2+1}$\\
\hline
Notation & $\hat{Z}_x$ & $\hat{Z}_y$ & $\hat{X}$ & $\hat{Y}$ & $\hat{I}$
\end{tabular}
%\caption{Finite dimensional representations of basic differential operators in two dimensions.}
\end{center}
\begin{center}
\small Table 2. Finite dimensional representations of basic differential operators in two dimensions.
\normalsize
\end{center}

Denote the right hand side of equation~\eqref{Ex2eq1} by $$f(x,y)=4(y^2-x^2)\sin (1-x^2-y^2)-2xy \cos (1-x^2-y^2).$$

Equation~\eqref{Ex2eq1} is represented by the system of linear equations
\begin{equation}
\label{eq:KvF}
\hat{K} v_{\sigma}=\hat{F},
\end{equation}
where 
\begin{equation}
\hat{K} =(\hat{I}-\hat{X}^2-\hat{Y}^2)(\hat{Z}_x^2-\hat{Z}_y^2+\hat{Y}\hat{Z}_x)-4\hat{X}\hat{Z}_x+4\hat{Y}\hat{Z}_y-2\hat{X}\hat{Y},
\end{equation}
\begin{equation*}
\hat{F} =(f(x_0,y_0), f(x_1,y_0), \ldots, f(x_{n_1},y_0), f(x_0, y_1), f(x_1,y_1),\ldots, f(x_{n_1},y_{n_2}))^T,
\end{equation*}
and the unknown vector $v_\sigma \in \R^{(n_1+1)(n_2+1)}$.

After solving system of linear equations~\eqref{eq:KvF}, we compute the vector $u_\sigma$ of approximate values of $u$ as follows:
$$
u_\sigma=(\hat{I}-\hat{X}^2-\hat{Y}^2)v_\sigma.
$$

We have compared the approximate solutions with the exact solution by computing the maximal error $E(n_1,n_2)=\max_{(i,j)\in \mathcal{I}} |[u_\sigma]_{\ast(i,j)}-u(x_i,y_j)|$ and the average error $E_{a}(n_1,n_2)=\frac{1}{(n_1+1)(n_2+1)}\sum_{(i,j)\in \mathcal{I}} |[u_\sigma]_{\ast(i,j)}-u(x_i,y_j)|$. The values of $E$ and $E_a$ for different number of partition points are given in table~3, and the plot of the approximate solution $u_\sigma$ is given in figure~1.

\pagebreak

\begin{table}[h]
\begin{center}
\begin{tabular}{c|c|c}
\label{table3}
$(n_1, n_2)$ & $(10,10)$ & $(15,15)$\\ 
\hline
$E(n_1,n_2)$ & 0.0064 & 0.002\\
$E_a(n_1,n_2)$ & 2.56e-04 & 2.63e-05
\end{tabular}
%\caption{Maximal and average error of the Calogero method for problem~\eqref{eq:lala},~\eqref{eq:lalaBC}.}
\end{center}
\end{table}
\vspace{-25pt}
\begin{center}
\small Table 3. Maximal and average error of the Lie-algebraic method for problem~\eqref{eq:lala},~\eqref{eq:lalaBC}.
\normalsize
\end{center}

\begin{figure}[h]
\begin{center} 
\includegraphics[width=10cm]{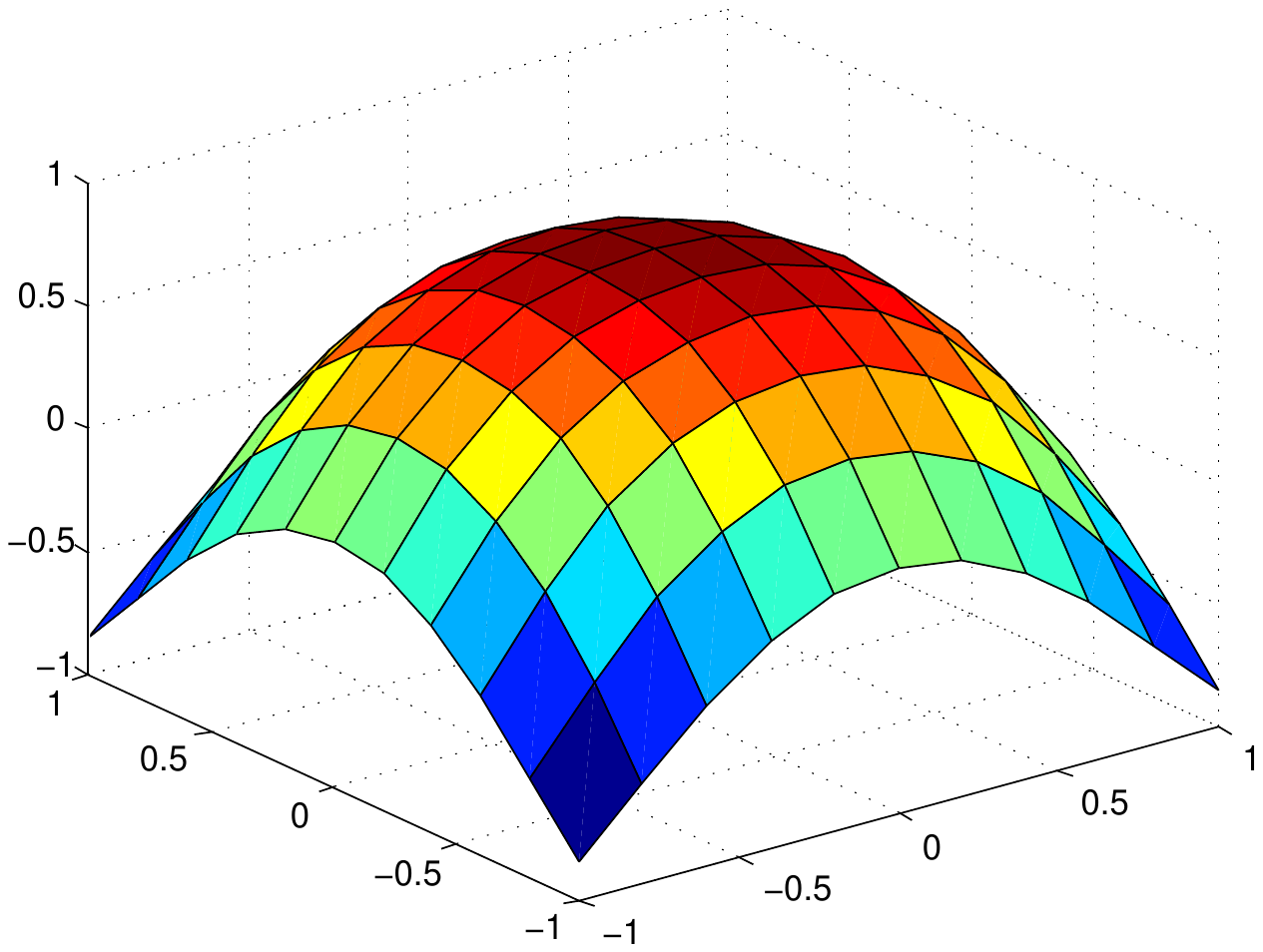}
\label{ppict1}
%\caption{Plot of the approximate solution $u_{\sigma}$ to problem ~\eqref{eq:lala},~\eqref{eq:lalaBC}.}
\end{center}
\end{figure}
\begin{center}
\small Figure~1. Plot of the approximate solution $u_{\sigma}$ to problem ~\eqref{eq:lala},~\eqref{eq:lalaBC}.
\normalsize
\end{center}
\end{example}

 {\small
%\begin{thebibliography}{99}

%\end{thebibliography}

\vskip6pt
${}^{\ast}${Concordia College}

{901 8th Str. South, Moorhead, MN 56562, USA}

\vskip20pt

${}^{\ast\ast}${Ivan Franko National University of Lviv}

{Universytetska Str. 1, 79000 Lviv, Ukraine}

\end{document}